# MULTIVARIABLE APPROXIMATE CARLEMAN-TYPE THEOREMS FOR COMPLEX MEASURES[1]

BY ISABELLE CHALENDAR AND JONATHAN R. PARTINGTON

*Université Lyon 1 and University of Leeds*

We prove a multivariable approximate Carleman theorem on the determination of complex measures on $\mathbb{R}^n$ and $\mathbb{R}^n_+$ by their moments. This is achieved by means of a multivariable Denjoy–Carleman maximum principle for quasi-analytic functions of several variables. As an application, we obtain a discrete Phragmén–Lindelöf-type theorem for analytic functions on $\mathbb{C}^n_+$.

**1. Introduction.** The main issue discussed in this paper is the determination of a complex measure on $\mathbb{R}^n$ from its moments. This is strongly linked to the theory of quasi-analytic classes, in which the partial derivatives of a function determine the function uniquely. Some work in this direction has recently been presented in [3, 4]. One of the aims of this paper is to present a generalization of some results from [3], describing the consequences when the moments of two measures differ by a set of values that does not grow too quickly. Moreover, our results are valid for complex measures, not just the positive measures discussed in the earlier work. Analogous results are presented for measures supported on the positive cone $\mathbb{R}^n_+$ (all such notation is defined below). An application of this is given in the form of a discrete Phragmén–Lindelöf-type theorem, generalizing some results of [8] which apply in the one-dimensional case. This paper also contains a multivariable generalization of [2], where a one-dimensional Denjoy–Carleman maximum principle as well as an approximate Carleman theorem are proved.

We adopt the conventions $\mathbb{N} = \{0, 1, 2, \ldots\}$, $\mathbb{R}_+ = \{x \in \mathbb{R} : x > 0\}$, $\mathbb{R}_- = \{x \in \mathbb{R} : x < 0\}$ and $\mathbb{C}_+ = \{z \in \mathbb{C} : \operatorname{Re} z > 0\}$. Let $n$ be a positive integer. We supply $\mathbb{R}^n$ and $\mathbb{C}^n$ with the standard inner product $\langle \cdot, \cdot \rangle$ and corresponding

Received February 2005; revised August 2005.
[1]Supported by the EPSRC.

norm $\|\cdot\|$. We denote by $\mathcal{M}(\mathbb{R}^n)$ the set of all positive Borel measures $\mu$ on $\mathbb{R}^n$ such that

$$\int_{\mathbb{R}^n} \|x\|^d \, d\mu(x) < \infty \qquad \text{for all } d \geq 0$$

and let $\mathcal{M}_c(\mathbb{R}^n)$ be the set of all complex Borel measures $\mu$ such that $|\mu| \in \mathcal{M}(\mathbb{R}^n)$. For $\alpha = (\alpha_1, \ldots, \alpha_n) \in \mathbb{N}^n$, set $|\alpha| = \alpha_1 + \cdots + \alpha_n$ and $\alpha! = \alpha_1! \cdots \alpha_n!$. For $f : \mathbb{R}^n \to \mathbb{C}$, $f^{(\alpha)}$ denotes $\frac{\partial^{\alpha_1}}{\partial x_1^{\alpha_1}} \cdots \frac{\partial^{\alpha_n}}{\partial x_n^{\alpha_n}} f$. As usual, for $z = (z_1, \ldots, z_n) \in \mathbb{C}^n$, $z^\alpha$ denotes $z_1^{\alpha_1} \cdots z_n^{\alpha_n}$.

Given a product $I = I_1 \times \cdots \times I_n$ of subintervals of $\mathbb{R}$ (bounded or unbounded) containing 0, and an $n$-tuple $M = (M(1,k))_{k \geq 0}, \ldots, (M(n,k))_{k \geq 0}$ of sequences of positive numbers, write $\mathcal{C}_I(M)$ for the family of all $C^\infty$-functions $f : I \to \mathbb{C}$ satisfying

$$(1) \quad |f^{(\alpha)}(x)| \leq c_f \rho_f^{|\alpha|} \prod_{j=1}^{n} M(j, \alpha_j), \qquad [x \in I, \alpha = (\alpha_1, \ldots, \alpha_n) \in \mathbb{N}^n],$$

where $c_f$ and $\rho_f$ are constants depending on $f$.

In the sequel, we will say that a sequence $(m_k)_{k \geq 0}$ of positive reals is a *Carleman sequence* if

$$m_0 = 1, \qquad m_k^2 \leq m_{k-1} m_{k+1} \quad \text{and} \quad \sum_{k \geq 1} m_k^{-1/k} = \infty.$$

The structure of this paper is as follows. In Section 2, we derive a multivariable Denjoy–Carleman maximum principle (Theorem 2.3) by means of a multivariable version of Bernstein's inequality for functions of exponential type (Theorem 2.2). This provides an extension of the results in [2].

As an application of the ideas of Section 2, we obtain in Section 3 a multivariable approximate Carleman theorem on $\mathbb{R}^n$ (Theorem 3.2), which holds even for complex measures. The methods used include a Paley–Wiener–Schwartz-type theorem for Fourier transforms of distributions.

In Section 4, we derive analogous results for complex measures supported on $\mathbb{R}_+^n$ (see Theorem 4.1). Even in one dimension, the results provide sharper forms of Theorem 2.1 of [8]; we then use the methods of [8] to derive a multivariable discrete Phragmén–Lindelöf-type theorem which, in the one-dimensional case, extends Theorem 4.1 of the same paper.

**2. A multivariable Denjoy–Carleman maximum principle.** We recall the multivariable Denjoy–Carleman theorem:

THEOREM 2.1 ([4], page 155, and [6]). *Let $n$ be a positive integer. For $j = 1, \ldots, n$, set $I_j = [-R_j, R_j]$, where $R_j > 0$ and $I = I_1 \times \cdots \times I_n$. Let $f \in \mathcal{C}_I(M)$ where $M$ is an $n$-tuple of Carleman sequences. If $f^{(\alpha)}(0) = 0$ for all $\alpha \in \mathbb{N}^n$, then $f$ is identically equal to $0$ on $I$.*

In fact, Theorem 2.1 is a slightly different version of the result given in [4], but the proof implies the version stated here.

We say that an entire function $g: \mathbb{C}^n \to \mathbb{C}$ is *of exponential type at most* $\tau = (\tau_1, \ldots, \tau_n)$, where $\tau_k > 0$ for all $k$, if for all $\varepsilon > 0$, there is a constant $A_\varepsilon$ such that

$$|g(z)| \leq A_\varepsilon e^{(\tau_1+\varepsilon)|z_1|+\cdots+(\tau_n+\varepsilon)|z_n|} \qquad \text{for each } z = (z_1, \ldots, z_n) \in \mathbb{C}^n.$$

The following multivariable version of Bernstein's theorem is easily deduced from the single-variable case:

THEOREM 2.2. *Let $g: \mathbb{C}^n \to \mathbb{C}$ be a function of exponential type at most $\tau = (\tau_1, \ldots, \tau_n)$. Then for each $k$ with $1 \leq k \leq n$, the function $\frac{\partial g}{\partial z_k}$ is also of exponential type at most $\tau$. Further, if $g$ is bounded on $\mathbb{R}^n$, then so is $\frac{\partial g}{\partial z_k}$, and*

$$\sup_{x \in \mathbb{R}^n} \left| \frac{\partial g}{\partial z_k}(x) \right| \leq \tau_k \sup_{x \in \mathbb{R}^n} |g(x)|.$$

PROOF. Without loss of generality, we take $k = n$. For each $(z_1, \ldots, z_{n-1}) \in \mathbb{C}^{n-1}$, let $g_{(z_1,\ldots,z_{n-1})}: \mathbb{C} \to \mathbb{C}$ be defined by $g_{(z_1,\ldots,z_{n-1})}(z_n) = g(z_1, \ldots, z_n)$. Then for each $\varepsilon > 0$, there is an $A_\varepsilon > 0$ such that

$$|g_{(z_1,\ldots,z_{n-1})}(z_n)| \leq A_\varepsilon e^{(\tau_1+\varepsilon)|z_1|+\cdots+(\tau_{n-1}+\varepsilon)|z_{n-1}|} e^{(\tau_n+\varepsilon)|z_n|}.$$

Using the Cauchy integral formula for $g'_{(z_1,\ldots,z_{n-1})}$ with the circle centered at $z_n$ and of radius 1, we obtain

$$|g'_{(z_1,\ldots,z_{n-1})}(z_n)| \leq A_\varepsilon e^{(\tau_1+\varepsilon)|z_1|+\cdots+(\tau_{n-1}+\varepsilon)|z_{n-1}|} e^{(\tau_n+\varepsilon)(|z_n|+1)}$$

$$= A_\varepsilon e^{\tau_n+\varepsilon} e^{(\tau_1+\varepsilon)|z_1|+\cdots+(\tau_n+\varepsilon)|z_n|}.$$

Clearly this implies that $\frac{\partial g}{\partial z_n}$ is of exponential type at most $\tau$. By Bernstein's theorem ([1], Theorems 2.4.1 and 11.1.2) applied to $g_{z_1,\ldots,z_{n-1}}$, for $z_1, \ldots, z_{n-1}$ real, we obtain

$$\sup_{x \in \mathbb{R}} \left| \frac{\partial g}{\partial z_n}(z_1, \ldots, z_{n-1}, x) \right| \leq \tau_n \sup_{x \in \mathbb{R}} |g(z_1, \ldots, z_{n-1}, x)|,$$

from which the result follows. □

THEOREM 2.3. *Let $f \in \mathcal{C}_I(M)$, where $I = \mathbb{R}^n$ and $M$ is an $n$-tuple of Carleman sequences. Suppose that there exist $C_1, \ldots, C_n > 0$ such that for each $\varepsilon > 0$, there is a constant $A_\varepsilon$ such that*

$$|f^{(\alpha)}(0)| \leq A_\varepsilon (C_1 + \varepsilon)^{\alpha_1} \cdots (C_n + \varepsilon)^{\alpha_n} \qquad \text{for all } \alpha \in \mathbb{N}^n.$$

*Then*

$$\sup_{x \in \mathbb{R}^n} |f^{(\alpha+\beta)}(x)| \leq C_1^{\beta_1} \cdots C_n^{\beta_n} \sup_{x \in \mathbb{R}^n} |f^{(\alpha)}(x)| \qquad \text{for all } \alpha, \beta \in \mathbb{N}^n.$$

PROOF. Define $h: \mathbb{C}^n \to \mathbb{C}$ by

$$h(z) = \sum_{\alpha \in \mathbb{N}^n} \frac{1}{\alpha!} f^{(\alpha)}(0) z^\alpha.$$

Thus for each $\varepsilon > 0$, we have

$$|h(z)| \leq \sum_{\alpha \in \mathbb{N}^n} A_\varepsilon (C_1 + \varepsilon)^{\alpha_1} \cdots (C_n + \varepsilon)^{\alpha_n} \frac{1}{\alpha!} |z_1|^{\alpha_1} \cdots |z_n|^{\alpha_n}$$

$$= A_\varepsilon e^{(C_1+\varepsilon)|z_1| + \cdots + (C_n+\varepsilon)|z_n|}.$$

That is, $h$ is an entire function of exponential type at most $(C_1, \ldots, C_n)$. For $\beta \in \mathbb{N}^n$, we have

$$h^{(\beta)}(z) = \sum_{\alpha \in \mathbb{N}^n} \frac{f^{(\alpha)}(0)}{\alpha!} \frac{\partial^\beta}{\partial z^\beta}(z^\alpha) = \sum_{\gamma \in \mathbb{N}^n} f^{(\beta+\gamma)}(0) \frac{z^\gamma}{\gamma!},$$

writing $\alpha = \beta + \gamma$. Therefore,

$$|h^{(\beta)}(z)| \leq A_\varepsilon \sum_{\gamma \in \mathbb{N}^n} (C_1 + \varepsilon)^{\beta_1+\gamma_1} \cdots (C_n + \varepsilon)^{\beta_n+\gamma_n} |z_1|^{\gamma_1} \cdots |z_n|^{\gamma_n} \frac{1}{\gamma!}$$

$$\leq A_\varepsilon (C_1 + \varepsilon)^{\beta_1} \cdots (C_n + \varepsilon)^{\beta_n} e^{(C_1+\varepsilon)|z_1| + \cdots + (C_n+\varepsilon)|z_n|}.$$

So given $R > 0$, we have

$$\sup_{x \in [-R,R]^n} |h^{(\beta)}(x)| \leq A_\varepsilon (C_1 + \varepsilon)^{\beta_1} \cdots (C_n + \varepsilon)^{\beta_n} e^{R((C_1+\varepsilon) + \cdots + (C_n+\varepsilon))}.$$

Note that since $(M(j, k+1)/M(j, k))_{k \geq 0}$ is an increasing sequence for each $j$, we have $M(j, 1)^{\beta_j} \leq M(j, \beta_j)$. Hence,

$$\sup_{x \in [-R,R]^n} |h^{(\beta)}(x)| \leq A_\varepsilon \prod_{j=1}^n \frac{(C_j + \varepsilon)^{\beta_j}}{M(j,1)^{\beta_j}} e^{R((C_1+\varepsilon) + \cdots + (C_n+\varepsilon))} \prod_{j=1}^n M(j, \beta_j)$$

$$\leq A_\varepsilon e^{R((C_1+\varepsilon) + \cdots + (C_n+\varepsilon))} \rho^{|\beta|} \prod_{j=1}^n M(j, \beta_j),$$

where $\rho = \max\{(C_j + \varepsilon)/(M(j,1)) : 1 \leq j \leq n\}$. Therefore, $h \in \mathcal{C}_{[-R,R]^n}(M)$ and so $f - h \in \mathcal{C}_{[-R,R]^n}(M)$, with $(f - h)^{(\alpha)}(0) = 0$ for each $\alpha \in \mathbb{N}^n$. By Theorem 2.1, $f - h$ vanishes identically on $[-R, R]^n$. This holds for all $R > 0$ and so $f \equiv h$ on $\mathbb{R}^n$.

By repeatedly applying Theorem 2.3, we obtain that $h^{(\alpha)}$ is of exponential type at most $(C_1, \ldots, C_n)$ and that

$$\sup_{x \in \mathbb{R}^n} |h^{(\alpha+\beta)}(x)| \leq C_1^{\beta_1} \cdots C_n^{\beta_n} \sup_{x \in \mathbb{R}^n} |h^{(\alpha)}(x)|.$$

Since $h_{|\mathbb{R}^n} = f$, the result follows. □

The following corollary is immediate and shows that Theorem 2.3 is an extension of Theorem 2.1:

COROLLARY 2.4. *Let $f$ satisfy the hypotheses of Theorem 2.3 and suppose that $\limsup_{|\alpha| \to \infty} |f^{(\alpha)}(0)|^{1/|\alpha|} = 0$. Then $f$ is constant.*

**3. An extension of the multivariable Carleman theorem on $\mathbb{R}^n$.** The following multivariable Carleman theorem is proved by De Jeu in [3]:

THEOREM 3.1. *Suppose that $\mu_1, \mu_2 \in \mathcal{M}(\mathbb{R}^n)$ satisfy*

$$(2) \qquad s(\alpha) := \int_{\mathbb{R}^n} x^\alpha \, d\mu_1(x) = \int_{\mathbb{R}^n} x^\alpha \, d\mu_2(x) \qquad \text{for all } \alpha \in \mathbb{N}^n$$

*and that the conditions*

$$\sum_{m=1}^{\infty} s(2me_j)^{-1/(2m)} = \infty, \qquad j = 1, \ldots, n,$$

*hold, where $e_j$ is the $j$th canonical basis vector of $\mathbb{R}^n$. Then $\mu_1 = \mu_2$.*

As an application of the ideas of the previous section, we obtain the following multivariable approximate Carleman theorem. Note that it applies even to complex measures and, in the case $n = 1$, it provides a strong generalization of [9].

THEOREM 3.2. *Let $C_1, \ldots, C_n > 0$. Suppose that $\mu_1, \mu_2 \in \mathcal{M}_c(\mathbb{R}^n)$ satisfy*

$$\int_{\mathbb{R}^n} x^\alpha \, d\mu_1(x) = \int_{\mathbb{R}^n} x^\alpha \, d\mu_2(x) + c(\alpha) \qquad \text{for all } \alpha \in \mathbb{N}^n,$$

*where for all $\varepsilon > 0$, there exists a constant $A_\varepsilon > 0$ such that*

$$|c(\alpha)| \leq A_\varepsilon (C_1 + \varepsilon)^{\alpha_1} \cdots (C_n + \varepsilon)^{\alpha_n}$$

*holds for all $\alpha \in \mathbb{N}^n$. Let*

$$s(\alpha) = \int_{\mathbb{R}^n} x^\alpha \, d(|\mu_1| + |\mu_2|)(x).$$

*Suppose that the conditions*

$$(3) \qquad \sum_{m=1}^{\infty} s(2me_j)^{-1/(2m)} = \infty, \qquad j = 1, \ldots, n,$$

*hold. Then $\mu_1 = \mu_2 + \sigma$, where $\sigma$ is a complex measure on $\mathbb{R}^n$ supported on $[-C_1, C_1] \times \cdots \times [-C_n, C_n]$.*

We immediately deduce the following strengthened version of Theorem 3.1:

COROLLARY 3.3. *Let $\mu_1$ and $\mu_2$ satisfy the hypotheses of Theorem 3.2, with $c(\alpha)$ satisfying $\limsup_{|\alpha|\to\infty} |c(\alpha)|^{1/|\alpha|} = 0$. Then $\mu_1 = \mu_2 + a\delta_0$ for some $a \in \mathbb{C}$. If, in addition, $\mu_1$ and $\mu_2$ are probability measures, then $\mu_1 = \mu_2$.*

The proof of Theorem 3.2 will require two preliminary results, which we now present. The one-dimensional version of the first result can be found in [7] and its extension follows by an inductive argument.

THEOREM 3.4. *Suppose that $h:\mathbb{C}^n \to \mathbb{C}$ is of exponential type at most $\tau = (\tau_1, \ldots, \tau_n)$ and that $|h|$ is bounded on $\mathbb{R}^n$ by a constant $D \geq 0$. Then*

$$(4) \quad |h(z)| \leq De^{\tau_1|\operatorname{Im} z_1| + \cdots + \tau_n|\operatorname{Im} z_n|} \qquad \text{for all } z = (z_1, \ldots, z_n) \in \mathbb{C}^n.$$

The next result is a version of the Paley–Wiener–Schwartz theorem and is a special case of Theorem 7.3.1 on page 181 of [5].

THEOREM 3.5. *Let $\tau_1, \ldots, \tau_n > 0$. If an entire function $h:\mathbb{C}^n \to \mathbb{C}$ satisfies (4), then it is the Fourier transform of a distribution supported on $[-\tau_1, \tau_1] \times \cdots \times [-\tau_n, \tau_n]$.*

We are now ready for the proof of Theorem 3.2.

PROOF OF THEOREM 3.2. For all $m \geq 0$ and $j = 1, \ldots, n$, define

$$M(j,m) = \frac{1}{m_0} \int_{\mathbb{R}^n} |x_j|^m \, d(|\mu_1| + |\mu_2|)(x),$$

where $m_0 = (|\mu_1| + |\mu_2|)(\mathbb{R}^n)$. Clearly, $M(j,0) = 1$, and on applying Hölder's inequality, we obtain

$$M(j,m)^2 \leq M(j,m-1)M(j,m+1) \qquad \text{for all } m \geq 1.$$

Moreover, note that

$$M(j,m) = \frac{1}{m_0} \int_{\mathbb{R}} |x_j|^m \, d\nu(x_j),$$

where

$$d\nu(x_j) = \int_{(x_1,\ldots,x_{j-1},x_{j+1},\ldots,x_n)\in\mathbb{R}^{n-1}} d(|\mu_1| + |\mu_2|)(x_1,\ldots,x_n).$$

Hence, using the calculations from [2], page 93, it follows that

$$\sum_{m=1}^{\infty} M(j,m)^{-1/m} = \infty \qquad \text{for each } j = 1, \ldots, n.$$

In other words, $M := (M(1,k)_{k\geq 0}, \ldots, M(n,k)_{k\geq 0})$ is an $n$-tuple of Carleman sequences.

Define $f : \mathbb{R}^n \to \mathbb{C}$ by

$$f(w) = \frac{1}{m_0} \int_{\mathbb{R}^n} e^{-i\langle w, x\rangle} \, d(\mu_1 - \mu_2)(x),$$

so that

$$f^{(\alpha)}(w) = \frac{1}{m_0} \int_{\mathbb{R}^n} (-ix_1)^{\alpha_1} \cdots (-ix_n)^{\alpha_n} e^{-i\langle w, x\rangle} \, d(\mu_1 - \mu_2)(x),$$

where $\alpha = (\alpha_1, \ldots, \alpha_n)$ and $x = (x_1, \ldots, x_n)$. Therefore, we have

$$|f^{(\alpha)}(w)| \leq \frac{1}{m_0} \int_{\mathbb{R}^n} |x_1|^{\alpha_1} \cdots |x_n|^{\alpha_n} \, d(|\mu_1| + |\mu_2|)(x).$$

Using the generalized Hölder inequality, we obtain

$$|f^{(\alpha)}(w)| \leq M(1, n\alpha_1)^{1/n} \cdots M(1, n\alpha_n)^{1/n}.$$

Now, for each $j = 1, \ldots, n$, set $\widetilde{M}(j, \alpha_j) = M(j, n\alpha_j)^{1/n}$. It is clear that $\widetilde{M}(j, 0) = 1$. Moreover, for all $k \geq 0$, $\widetilde{M}(j, k)^2 \leq \widetilde{M}(j, k-1)\widetilde{M}(j, k+1)$ since $(\frac{M(j, m+1)}{M(j, m)})_{m \geq 0}$ is increasing. The fact that $\sum_{m=1}^\infty \widetilde{M}(j, m)^{-1/m} = \infty$ follows from Lemma 2.2 of [3] and the fact that $(M(j, m)^{1/m})_{m \geq 1}$ is increasing. At this stage, we have proved that $\widetilde{M} := ((\widetilde{M}(1, k))_{k \geq 0}, \ldots, (\widetilde{M}(n, k))_{k \geq 0})$ is an $n$-tuple of Carleman sequences and that $f \in \mathcal{C}_{\mathbb{R}^n}(\widetilde{M})$.

Now, note that $f^{(\alpha)}(0) = \frac{(-i)^{|\alpha|}}{m_0} c(\alpha)$. Therefore, for all $\varepsilon > 0$, there exists $A_\varepsilon > 0$ such that

$$|f^{(\alpha)}(0)| \leq A_\varepsilon (C_1 + \varepsilon)^{\alpha_1} \cdots (C_n + \varepsilon)^{\alpha_n}.$$

As in the proof of Theorem 2.3, there is an entire function $h : \mathbb{C}^n \to \mathbb{C}$ of exponential type at most $(C_1, \ldots, C_n)$ such that $h_{|\mathbb{R}^n} = f$ and $|h(z)| \leq 1$ for all $z \in \mathbb{R}^n$. By Theorem 3.4, we deduce

$$|h(z)| \leq e^{C_1 |\operatorname{Im}(z_1)| + \cdots + C_n |\operatorname{Im}(z_n)|} \qquad \text{for all } z = (z_1, \ldots, z_n) \in \mathbb{C}^n.$$

Hence, using Theorem 3.5, $h$ is the Fourier–Laplace transform of a distribution $u$ supported on $[-C_1, C_1] \times \cdots \times [-C_n, C_n]$. Thus, $f$ is just the Fourier transform of $u$. But $f$ was defined as the Fourier transform of $\frac{\mu_1 - \mu_2}{m_0}$. So, by the uniqueness theorem for Fourier transforms of tempered distributions on $\mathbb{R}^n$, $u = \frac{\mu_1 - \mu_2}{m_0}$. In particular, $\mu_1 - \mu_2$ is supported on $[-C_1, C_1] \times \cdots \times [-C_n, C_n]$, as required. $\square$

Here is a corollary involving a weaker condition than (2) but valid only under a restrictive condition on the measures $\mu_1$ and $\mu_2$:

COROLLARY 3.6. *Let $C_1, \ldots, C_n > 0$. Suppose that $\mu_1, \mu_2 \in \mathcal{M}(\mathbb{R}^n)$ satisfy $\mu_1 \geq \mu_2$ and that*

$$(5) \quad s(j, m) := \int_{\mathbb{R}^n} x_j^m \, d\mu_1(x) = \int_{\mathbb{R}^n} x_j^m \, d\mu_2(x) + c(j, m) \qquad \text{for all } m \in \mathbb{N},$$

*where for all $\varepsilon > 0$, there exists a constant $A_\varepsilon > 0$ such that*

$$|c(j, m)| \leq A_\varepsilon (C_j + \varepsilon)^{\alpha m}$$

*for all $m \in \mathbb{N}$. Suppose further that the conditions*

$$\sum_{m=1}^{\infty} s(j, 2m)^{-1/(2m)} = \infty, \qquad j = 1, \ldots, n,$$

*hold. Then $\mu_1 = \mu_2 + \sigma$, where $\sigma$ is a positive measure on $\mathbb{R}^n$ supported on $[-C_1, C_1] \times \cdots \times [-C_n, C_n]$.*

PROOF. Following the lines of the proof of Theorem 3.2, since $\mu_1 - \mu_2 \geq 0$, we first obtain

$$|f^{(\alpha)}(w)| \leq \frac{1}{m_0} \int_{\mathbb{R}^n} |x_1|^{\alpha_1} \cdots |x_n|^{\alpha_n} \, d(\mu_1 - \mu_2)(x)$$

and then, using the generalized Hölder inequality, it follows that

$$|f^{(\alpha)}(w)| \leq M(1, n\alpha_1)^{1/n} \cdots M(1, n\alpha_n)^{1/n}.$$

The proof ends in the same way as the proof of Theorem 3.2. □

REMARK 3.1. If in Theorem 3.2 we take $\mu_1$ and $\mu_2$ to be in $\mathcal{M}(\mathbb{R}^n)$ rather than in $\mathcal{M}_c(\mathbb{R}^n)$, then, clearly, we may replace condition (3) by

$$\sum_{m=1}^{\infty} \tilde{s}(2me_j)^{-1/(2m)} = \infty, \qquad j = 1, \ldots, n,$$

where

$$\tilde{s}(\alpha) = \int_{\mathbb{R}^n} x^\alpha \, d\mu_1(x).$$

## 4. An approximate multivariable Carleman theorem on $\mathbb{R}_+^n$ with applications.

4.1. *An extension of the multivariable Carleman theorem on $\mathbb{R}_+^n$.*

THEOREM 4.1. *Let $C_1, \ldots, C_n > 0$. Suppose that $\mu_1, \mu_2 \in \mathcal{M}_c(\mathbb{R}_+^n)$ satisfy*

$$\int_{\mathbb{R}_+^n} x^\alpha \, d\mu_1(x) = \int_{\mathbb{R}_+^n} x^\alpha \, d\mu_2(x) + c(\alpha) \qquad \text{for all } \alpha \in \mathbb{N}^n,$$

*where for all $\varepsilon > 0$, there exists a constant $A_\varepsilon > 0$ such that*

$$|c(\alpha)| \leq A_\varepsilon (C_1 + \varepsilon)^{\alpha_1} \cdots (C_n + \varepsilon)^{\alpha_n}$$

*holds for all $\alpha \in \mathbb{N}^n$. Let*

$$s(\alpha) = \int_{\mathbb{R}_+^n} x^\alpha \, d(|\mu_1| + |\mu_2|)(x).$$

*Suppose further that the conditions*

(6) $$\sum_{m=1}^\infty s(me_j)^{-1/(2m)} = \infty, \qquad j = 1, \ldots, n,$$

*hold. Then $\mu_1 = \mu_2 + \sigma$, where $\sigma$ is a complex measure on $\mathbb{R}_+^n$ supported on $[0, C_1] \times \cdots \times [0, C_n]$.*

PROOF. We define measures $\nu_1$ and $\nu_2$ on $\mathbb{R}^n$ by

$$d\nu_k(t_1, \ldots, t_n) = d\mu_k(t_1^2, \ldots, t_n^2) \qquad \text{for } k = 1, 2.$$

If $\alpha \in \mathbb{N}^n$ satisfies $\alpha = 2\beta$ for some $\beta \in \mathbb{N}^n$, then

$$s'(\alpha) := \int_{\mathbb{R}^n} t^\alpha d(|\nu_1| + |\nu_2|)(t) = 2 \int_{\mathbb{R}_+^n} x^\beta d(|\mu_1| + |\mu_2|)(x),$$

otherwise $s'(\alpha) = 0$. We write

$$\int_{\mathbb{R}^n} t^\alpha \, d\nu_1(t) = \int_{\mathbb{R}^n} t^\alpha \, d\nu_2(t) + c'(\alpha).$$

Note that

$$|c'(\alpha)| \leq A_\varepsilon (C_1 + \varepsilon)^{\alpha_1/2} \cdots (C_n + \varepsilon)^{\alpha_n/2}$$

for all $\alpha \in \mathbb{N}^n$. Also

$$\sum_{m=1}^\infty s'(2me_j)^{-1/(2m)} = \sum_{m=1}^\infty (2s(me_j))^{-1/(2m)} = \infty, \qquad j = 1, \ldots, n.$$

Hence, by Theorem 3.2, $\nu_1 = \nu_2 + \sigma'$, where $\sigma'$ is a complex measure supported on $[-\sqrt{C_1}, \sqrt{C_1}] \times \cdots \times [-\sqrt{C_n}, \sqrt{C_n}]$. It follows that $\mu_1 - \mu_2$ is supported on $[0, C_1] \times \cdots \times [0, C_n]$. $\square$

The following corollary is immediate, and provides the multivariable Carleman theorem for $\mathbb{R}_+^n$ proved in [3], Theorem. 5.1:

COROLLARY 4.2. *Let $\mu_1$ and $\mu_2$ satisfy the hypotheses of Theorem 4.1, with $c(\alpha)$ satisfying $\limsup_{|\alpha| \to \infty} |c(\alpha)|^{1/|\alpha|} = 0$. Then $\mu_1 = \mu_2 + a\delta_0$ for some $a \in \mathbb{C}$. If, in addition, $\mu_1$ and $\mu_2$ are probability measures, then $\mu_1 = \mu_2$.*

REMARK 4.1. If in Theorem 4.1 we take $\mu_1$ and $\mu_2$ to be in $\mathcal{M}(\mathbb{R}_+^n)$ rather than in $\mathcal{M}_c(\mathbb{R}_+^n)$, then, clearly, we may replace condition (6) by

$$\sum_{m=1}^{\infty} \tilde{s}(me_j)^{-1/(2m)} = \infty, \qquad j = 1, \ldots, n,$$

where

$$\tilde{s}(\alpha) = \int_{\mathbb{R}_+^n} x^\alpha \, d\mu_1(x).$$

4.2. *A discrete Phragmén–Lindelöf theorem.* In this section, we present an application of Theorem 4.1; in the case $n = 1$, we obtain a generalization of [8], Theorem 4.1, including a simplification of the original one-dimensional proof.

We begin with an $n$-dimensional Phragmén–Lindelöf theorem, a stronger form of which can be found in [10]. See also [11].

THEOREM 4.3 ([10], page 303). *Let $f : \overline{\mathbb{C}_+^n} \to \mathbb{C}$ be continuous and holomorphic on $\mathbb{C}_+^n$ and let $a > 0$. Suppose that $|f|$ is bounded by $M$ on the set $(i\mathbb{R})^n$ and that for all $\varepsilon > 0$, there is a constant $c_\varepsilon > 0$ such that*

$$|f(z)| \leq c_\varepsilon e^{(a+\varepsilon)\|z\|} \qquad \text{for all } z \in \overline{\mathbb{C}_+^n},$$

*where $\|\cdot\|$ denotes any norm on $\mathbb{C}^n$. Then*

$$|f(z)| \leq M e^{a\|\operatorname{Re} z\|} \qquad \text{for all } z \in \overline{\mathbb{C}_+^n}.$$

The following theorem is of interest even in the case $n = 1$ which we discuss separately later:

THEOREM 4.4. *Let $C_1, \ldots, C_n \geq 1$. Let $f : \overline{\mathbb{C}_+^n} \to \mathbb{C}$ be continuous and holomorphic on $\mathbb{C}_+^n$ and bounded on each of the sets*

$$E_k := \{z = (z_1, \ldots, z_n) \in \overline{\mathbb{C}_+^n} : 0 \leq \operatorname{Re} z_k \leq 1\} \qquad \text{for } k = 1, \ldots, n.$$

*For $\alpha \in \mathbb{N}^n$, define*

$$M(\alpha) = \sup\{|f(z)| : 0 \leq \operatorname{Re} z_1 \leq \alpha_1, \ldots, 0 \leq \operatorname{Re} z_n \leq \alpha_n\}.$$

*Suppose that $M(\alpha) < \infty$ for all $\alpha$ and that*

$$(7) \qquad \sum_{m=1}^{\infty} M(me_j)^{-1/(2m)} = \infty, \qquad j = 1, \ldots, n.$$

*If for every $\varepsilon > 0$, there is a constant $A_\varepsilon > 0$ such that*

$$|f(\alpha)| \leq A_\varepsilon (C_1 + \varepsilon)^{\alpha_1} \cdots (C_n + \varepsilon)^{\alpha_n},$$

*then*

$$|f(z)| \leq M(0) C_1^{\operatorname{Re} z_1} \cdots C_n^{\operatorname{Re} z_n} \qquad \text{for all } z \in \overline{\mathbb{C}_+^n}.$$

*In particular, if* $\lim_{|\alpha|\to\infty} |f(\alpha)|^{1/|\alpha|} \leq 1$, *then* $|f|$ *is bounded by* $M(0)$ *on* $\overline{\mathbb{C}_+^n}$.

PROOF. For $z = (z_1, \ldots, z_n) = (x_1 + iy_1, \ldots, x_n + iy_n)$ with $x_k, y_k \in \mathbb{R}$ ($k = 1, \ldots, n$) and for all $\alpha \in \mathbb{N}^n$, set

$$g(z) = \frac{f(z)}{(1+z_1)\cdots(1+z_n)}.$$

Let

$$B(\alpha) = \sup_{0 \leq x_1 \leq \alpha_1, \ldots, 0 \leq x_n \leq \alpha_n} \left( \int_{\mathbb{R}^n} |g(z)|^2 \, dy_1 \cdots dy_n \right)^{1/2}.$$

Then

$$\begin{aligned} B(\alpha) &\leq \sup_{0 \leq x_1 \leq \alpha_1, \ldots, 0 \leq x_n \leq \alpha_n} \left( \int_{\mathbb{R}^n} \frac{|f(z)|^2 dy_1 \cdots dy_n}{|1+z_1|^2 \cdots |1+z_n|^2} \right)^{1/2} \\ &= M(\alpha) \pi^{n/2}. \end{aligned} \qquad (8)$$

Let us write $\operatorname{Re} z = (x_1, \ldots, x_n)$ and $\operatorname{Im} z = (y_1, \ldots, y_n)$. It follows from (8) that there exists a function $\Phi \in L^2(\mathbb{R}^n)$ such that

$$g(i \operatorname{Im} z) = \int_{\mathbb{R}^n} \Phi(\sigma) e^{\langle i \operatorname{Im} z, \sigma \rangle} \, d\sigma.$$

Moreover, its analytic extension to $\mathbb{C}_+^n$ is given by

$$g(z) = \int_{\mathbb{R}^n} \Phi(\sigma) e^{\langle z, \sigma \rangle} \, d\sigma.$$

Applying the multivariable Plancherel theorem, we obtain

$$\int_{\mathbb{R}^n} |g(z)|^2 \, dy_1 \cdots dy_n = (2\pi)^n \int_{\mathbb{R}^n} |\Phi(\sigma)|^2 e^{2\langle \operatorname{Re} z, \sigma \rangle} \, d\sigma$$

and thus

$$B(\alpha)^2 \geq (2\pi)^n \int_{\mathbb{R}^n} |\Phi(\sigma)|^2 e^{2\langle \alpha, \sigma \rangle} \, d\sigma. \qquad (9)$$

Let

$$\Psi(t_1, \ldots, t_n) = \Phi(\log t_1, \ldots, \log t_n) \qquad \text{for } t = (t_1, \ldots, t_n) \in \mathbb{R}_+^n.$$

We have

$$g(z + (1, \ldots, 1)) = \int_{\mathbb{R}^n} \Phi(\sigma) e^{\langle z + (1, \ldots, 1), \sigma \rangle} \, d\sigma.$$

Using the change of variables $\sigma = (\sigma_1, \ldots, \sigma_n) = (\log t_1, \ldots, \log t_n)$, this becomes

$$g(z + (1,\ldots,1)) = \int_{\mathbb{R}^n_+} \Psi(t) t_1^{z_1+1} \cdots t_n^{z_n+1} \frac{dt}{t_1 \cdots t_n} = \int_{\mathbb{R}^n_+} \Psi(t) t_1^{z_1} \cdots t_n^{z_n} \, dt.$$

For $\alpha \in \mathbb{N}^n$, set

$$C(\alpha) = \int_{\mathbb{R}^n_+} t_1^{\alpha_1} \cdots t_n^{\alpha_n} |\Psi(t)| \, dt$$

and note that

$$C(\alpha) = \int_{\mathbb{R}^n} e^{\sum_{k=1}^n (\alpha_k+1)\sigma_k} |\Phi(\sigma)| \, d\sigma$$

$$= \int_{\mathbb{R}^n_-} e^{\langle \alpha, \sigma \rangle} |\Phi(\sigma)| e^{\sum_{k=1}^n \sigma_k} \, d\sigma$$

$$+ \int_{\mathbb{R}^n_+} e^{\sum_{k=1}^n (\alpha_k+2)\sigma_k} |\Phi(\sigma)| e^{-\sum_{k=1}^n \sigma_k} \, d\sigma$$

$$\leq \left( \int_{\mathbb{R}^n_-} e^{2\langle \alpha, \sigma \rangle} |\Phi(\sigma)|^2 \, d\sigma \right)^{1/2} \left( \int_{\mathbb{R}^n_-} e^{2\sum_{k=1}^n \sigma_k} \, d\sigma \right)^{1/2}$$

$$+ \left( \int_{\mathbb{R}^n_+} e^{2\sum_{k=1}^n (\alpha_k+2)\sigma_k} |\Phi(\sigma)|^2 \, d\sigma \right)^{1/2} \left( \int_{\mathbb{R}^n_+} e^{-2\sum_{k=1}^n \sigma_k} \, d\sigma \right)^{1/2}.$$

By (9) and (8) we obtain

$$C(\alpha) \leq B(\alpha)(2\pi)^{-n/2} 2^{-n/2} + B(\alpha + (2,\ldots,2))(2\pi)^{-n/2} 2^{-n/2}$$

$$\leq 2(4\pi)^{-n/2} B(\alpha + (2,\ldots,2)) \leq 2^{1-n} M(\alpha + (2,\ldots,2)),$$

and hence, using (7), we have

$$\sum_{m=1}^\infty C(me_j)^{-1/(2m)} = \infty, \qquad j = 1, \ldots, n.$$

Clearly, $|g(z)| \leq |f(z)|$ for all $z \in \overline{\mathbb{C}^n_+}$, and so for every $\varepsilon > 0$, there is a constant $A_\varepsilon > 0$ such that

$$|g(\alpha)| \leq |f(\alpha)| \leq M(\alpha) \leq A_\varepsilon (C_1 + \varepsilon)^{\alpha_1} \cdots (C_n + \varepsilon)^{\alpha_n}.$$

Applying Theorem 4.1, we see that $\Psi$ is supported on $[0, C_1] \times \cdots \times [0, C_n]$, and so $\Phi$ is supported on $\mathcal{D}_n := (-\infty, \log C_1] \times \cdots \times (-\infty, \log C_n]$. Hence, we have

$$g(z) = \int_{\mathcal{D}_n} \Phi(\sigma) e^{\langle z, \sigma \rangle} \, d\sigma.$$

Since $f(z) = (z_1 + 1) \cdots (z_n + 1) g(z)$, it follows from the Cauchy–Schwarz inequality and the fact that $f$ is bounded on each set $E_k$ that

$$|f(z)| \leq K(|z|+1)^n C_1^{\operatorname{Re} z_1} \cdots C_n^{\operatorname{Re} z_n}$$

for some $K > 0$; also, $\sup\{|f(z)| : z \in (i\mathbb{R})^n\} = M(0)$. The conclusion is now an immediate consequence of Theorem 4.3 on taking $a = 1$ and $\|z\| = \sum_{k=1}^n (\log C_k + \delta)|z_k|$ and letting $\delta > 0$ tend to zero. $\square$

The following result generalizes [8], Theorem 4.1, which treats the special case when $(f(m))_{m \geq 1}$ is bounded. It follows immediately from Theorem 4.4.

COROLLARY 4.5. *Let $f : \overline{\mathbb{C}_+} \to \mathbb{C}$ be continuous and holomorphic on $\mathbb{C}_+$. For $m \in \mathbb{N}$, define*

$$M(m) = \sup\{|f(z)| : 0 \leq \operatorname{Re} z \leq m\}.$$

*Suppose that $M(m) < \infty$ for all $m$ and that*

$$\sum_{m=1}^\infty M(m)^{-1/(2m)} = \infty.$$

*If $\limsup_{m \to \infty} |f(m)|^{1/m} \leq 1$, then $|f|$ is bounded by $M(0)$ on $\overline{\mathbb{C}_+}$.*

**Acknowledgments.** The authors are grateful to the EPSRC for financial support.

UFR DE MATHÉMATIQUES
UNIVERSITÉ LYON 1
69622 VILLEURBANNE CEDEX
FRANCE
E-MAIL: chalenda@igd.univ-lyon1.fr

SCHOOL OF MATHEMATICS
UNIVERSITY OF LEEDS
LEEDS LS2 9JT
UK
E-MAIL: j.r.partington@leeds.ac.uk